\documentclass[12pt,leqno]{article}
\usepackage{amsmath,amsthm}

\def\init{\setcounter{equation}{0}}
\newtheorem{theorem}{Theorem}[section]

\newcommand{\R}{{\bf R}}
\newcommand{\C}{{\bf C}}

\newtheorem{lemma}{Lemma}[section]

\newcommand{\e}{{\varepsilon}}

\title{Inverse  problems  for Schr\"{o}dinger equations with
Yang-Mills potentials in
domains with  obstacles and the Aharonov-Bohm effect.
\author{G.Eskin, \ \ \  Department of Mathematics, UCLA,\\ Los Angeles,
CA 90095-1555, USA. \ E-mail: eskin@math.ucla.edu}
}

\begin{document}
\maketitle
\begin{abstract}
We study the inverse boundary value problems for the Schr\"{o}dinger equations 
with Yang-Mills potentials in a bounded domain $\Omega_0\subset\R^n$
containing finite number of smooth obstacles $\Omega_j,1\leq j \leq r$.
We prove that the Dirichlet-to-Neumann opeartor on $\partial\Omega_0$
determines the gauge equivalence class of the Yang-Mills potentials.
We also prove that the metric tensor can be recovered up to
a diffeomorphism that is identity on $\partial\Omega_0$.
\end{abstract}

\section{Introduction.}
\label{section 1}
\init

Let $\Omega_0$ be a smooth bounded domain in $\R^n$,
diffeomorphic to a ball, $\ n\geq 2,$
containing $r$ smooth nonintersecting obstacles $\Omega_j,\ 1\leq j\leq r$.
Consider the Schr\"{o}dinger equation in 
$\Omega=\Omega_0\setminus(\cup_{j=1}^r\overline{\Omega_j})$ with Yang-Mills
potentials
\begin{equation}                                  \label{eq:1.1}
\sum_{j=1}^n\left(-i\frac{\partial}{\partial x_j}I_m+A_j(x)\right)^2u
+V(x)u-k^2u=0
\end{equation}
with the boundary conditions
\begin{equation}                                 \label{eq:1.2}
u\left|_{\partial\Omega_j}\right.=0,\ \ \ 1\leq j\leq r,
\end{equation}
\begin{equation}                                 \label{eq:1.3}
u\left|_{\partial\Omega_0}\right.=f(x'),
\end{equation}
where  
$A_j(x),V(x),u(x)$ are $m\times m$ matrices, $I_m$ is the identity matrix in 
$\C^m$.  Let $G(\Omega)$ be the gauge group of all smooth nonsingular 
matrices in $\overline{\Omega}$.  Potentials $A(x)=(A_1,...,A_n),V$
and $A'(x)=(A_1',...,A_n'),V'(x)$ are called gauge equivalent if there 
exists $g(x)\in G(\Omega)$ such that
\begin{equation}                                 \label{eq:1.4}
A'(x)=g^{-1}Ag-ig^{-1}(x)\frac{\partial g}{\partial x},\ V'=g^{-1}Vg.
\end{equation}
Let $\Lambda$ be the Dirichlet-to-Neumann (D-to-N) operator on $\partial\Omega_0$,
i.e.
\[
\Lambda f=
(\frac{\partial u}{\partial \nu}
+i(A\cdot\nu)u)\left|_{\partial\Omega_0}^{\ }\right.,
\]
where $\nu=(\nu_1,...,\nu_n)$ is the unit outward normal to $\partial\Omega_0$
and $u(x)$ is the solution of (\ref{eq:1.1}), (\ref{eq:1.2}), (\ref{eq:1.3})).
We assume that the Dirichlet problem  
(\ref{eq:1.1}), (\ref{eq:1.2}), (\ref{eq:1.3}))
has a unique solution.
We shall say that the D-to-N operators $\Lambda$ and $\Lambda'$ are
gauge equivalent if there exists $g_0\in G(\Omega)$ such that
\[
\Lambda'=g_{0,\partial\Omega_0}^{\ }\Lambda g_{0,\partial\Omega_0}^{-1},
\]
where $g_{0,\partial\Omega_0}$ is the restriction of $g_0$ to 
$\partial\Omega_0$.  We shall prove the following theorem:
\begin{theorem}                        \label{theo:1.1}
Suppose that D-to-N operators $\Lambda'$ and $\Lambda$ corresponding 
to potentials $(A',V')$ and $(A,V)$ respectively are gauge equivalent
for all $k\in (k_0-\delta_0,k_0+\delta_0)$, where $k_0>0,\ \delta_0>0$.
Then potentials $(A',V')$ and $(A,V)$ are gauge equivalent too.
\end{theorem}
If we replace $A', V'$ by 
$A^{(1)}=g_0^{-1}A'g_0-ig_0^{-1}\frac{\partial g_0}{\partial x},\ 
V^{(1)}=g_0^{-1}Vg_0$ then $\Lambda=\Lambda_1$ where $\Lambda_1$ is
the D-to-N operator corresponding to $(A^{(1)},V^{(1)})$.
The proof of Theorem \ref{theo:1.1} gives that if $\Lambda=\Lambda_1$
then $(A,V)$ and $(A^{(1)},V^{(1)})$ are gauge equivalent with a gauge
$g\in G(\Omega)$ such 
that $g|_{\partial\Omega_0}=I_m$.
We shall denote the subgroup of $G(\Omega)$  consisting 
of $g$ such that 
$g(x)|_{\partial\Omega_0}=I_m$ by $G_0(\Omega)$.  In the case when 
$\Omega_0$ contains no obstacles Theorem \ref{theo:1.1} was proven in
[E] for $n\geq 3$ and in [E3] for $n=2$.  Note that the result of [E] is
stronger since it requires that $\Lambda=\Lambda^{(1)}$ for one value of
$k$ only.  In the case $n=2$ the proof of Theorem \ref{theo:1.1} is
simpler than that in [E3] since it does not rely on the uniqueness of 
the inversion of the non-abelian Radon transform.

We shall prove Theorem \ref{theo:1.1} in two steps.  In \S 2 we shall
 prove that $(A,V)$ and $(A^{(1)},V^{(2)})$ are locally gauge 
equivalent using the reduction to the inverse problem for the hyperbolic
equations as in 
[B], [B1],  [KKL], [KL], [E1], 
 and in \S 3 we shall prove the global gauge equivalence 
using the results of \S 2
and of [E2].  
Following Yang and Wu (see [WY]) one can describe the gauge 
equivalence class of $A=(A_1,...,A_n)$.  Fix a point $x^{(0)}\in 
\partial\Omega_0$ and consider all closed paths $\gamma$ in $\Omega$ 
starting and ending at $x^{(0)}$.  
Let $x=\gamma(\tau),\ 0\leq \tau\leq \tau_0,$
be a parametric equation of $\gamma,\ \gamma(0)=\gamma(\tau_0)=x^{(0)}$.
Consider the Cauchy problem for the system
\[
\frac{\partial}{\partial\tau}c(\tau,\gamma)=\frac{d\gamma(\tau)}{d\tau}\cdot
A(\gamma(\tau))c(\tau,\gamma),\ c(0,\gamma)=I_m.
\]
By the definition the gauge phase factor $c(\gamma,A)$ is $c(\tau_0,\gamma)$.
Therefore $A$ defines a map of the group of paths to $GL(m,\C)$.
The image of this map is a subgroup of $GL(m,\C)$ which is called the
holonomy group of $A$ (see  [Va]).
It is easy to show (c.f. \S 3)  that $c(\gamma,A^{(1)})=
c(\gamma,A^{(2)})$ for all closed paths $\gamma$ iff $A^{(1)}$ and
$A^{(2)}$ are gauge equivalent in 
$\Omega$. 
As it was shown
by Aharonov and Bohm [AB] the presence of distinct gauge equivalent
classes of potentials  
 can be detected in an experiment and this
phenomenon is called the Aharonov-Bohm effect.
In \$ 4 we consider the recovery of the Riemannian metrics from
the D-to-N operator in domains with obstacles. 

\section{Inverse problem for the hyperbolic system.}
\label{section 2}
\init

Consider two hyperbolic system:
\begin{equation}                               \label{eq:2.1}
L^{(p)}u=\frac{\partial^2}{\partial t^2}u^{(p)}
+\sum_{j=1}^n(-i\frac{\partial}{\partial x_j}I_m+A_j^{(p)}(x))^2u^{(p)}
+V^{(p)}(x)u^{(p)}=0,\ \ p=1,2,
\end{equation}
in
$\Omega\times (0,T_0)$ with zero initial conditions
\begin{equation}                            \label{eq:2.2}
u^{(p)}(x,0)=u_t^{(p)}(x,0)=0
\end{equation}
and the Dirichlet boundary conditions
\begin{equation}                           \label{eq:2.3}
u^{(p)}\left|_{\partial\Omega_j\times(0,T_0)}\right.=0,\ \
1\leq j\leq r,\ \ 
u^{(p)}\left|_{\partial\Omega_0\times(0,T_0)}\right.=f(x',t),\ \
p=1,2.
\end{equation}
Here $\Omega=\Omega_0\setminus(\cup_{j=1}^r\overline{\Omega_j})$
is the same as in \S 1, $A_j^{(p)}(x),\ 1\leq j\leq n,\ V^{(p)}(x),
\\ u^{(p)}(x,t),\ p=1,2,$ are smooth $m\times m$ matrices.
  As in \S 1 introduce D-to-N
operators
$\Lambda^{(p)}f=(\frac{\partial}{\partial\nu} 
+i\sum_{j=1}^n A_j^{(p)}\cdot\nu_j)u\left|_{\partial\Omega_0\times (0,T_0)}
\right.,\ p=1,2$. 

Making the Fourier transform in $t$ one can show that 
the D-to-N operator for (\ref{eq:2.1}) when $T_0=\infty$
determines the D-to-N operator for (\ref{eq:1.1})
for all $k$ except a discrete set,  and vice versa.

We shall prove the following theorem:
\begin{theorem}                           \label{theo:2.1}
Suppose $\Lambda^{(1)}=\Lambda^{(2)}$ and 
$T_0>\max_{x\in\overline{\Omega}}d(x,\partial\Omega_0)$
where $d(x,\partial\Omega_0)$ is the distance in $\overline{\Omega}$
from $x\in \overline{\Omega}$ to $\partial\Omega_0$.
Then potentials $A_j^{(1)}(x),1\leq j\leq n,\ V^{1}(x)$ and
$A_j^{(2)}(x),1\leq j\leq n,\ V^{(2)}(x)$ are gauge equivalent
in $\overline{\Omega}$,  i.e.  (\ref{eq:1.4}) holds with
$g\in G_0(\overline{\Omega})$.
\end{theorem}

Note that Theorem \ref{theo:2.1}  implies Theorem \ref{theo:1.1}.
We can consider a more general than (\ref{eq:2.1}) equation when
the Eucleadian metric is replaced by an arbitrary Riemannian metric:

\begin{eqnarray}                                   \label{eq:2.4}
\frac{\partial^2 u^{(p)}}
{\partial t^2} +\sum_{j,k=1}^n\frac{1}{\sqrt{g_p(x)}}
(-i\frac{\partial}{\partial x_j}I_m+A_j^{(p)}(x))\sqrt{g_p(x)}g_p^{jk}(x)
(-i\frac{\partial}{\partial x_j}I_m
\nonumber
\\
+A_k^{(p)}(x))u^{(p)}+V^{(p)}(x)u^{(p)}(x,t)=0,
\end{eqnarray}
where $\|g_p^{jk}(x)\|^{-1}$ are metric tensors in 
$\overline{\Omega}^{(p)}, g_p(x)=\det\|g_p^{jk}\|^{-1},\ 
A_j^{(p)}(x),V^{(p)}(x)$ are the same as in (\ref{eq:2.1}),
$\Omega^{(p)}=\Omega_0\setminus\overline{\Omega}_p',
\ \Omega_p'=\cup_{j=1}^r\Omega_{jp}$.
Let $\Gamma$ be an open subset of $\partial\Omega_0$ and let
$0< T<T_0$ be small.  Denote by $\Delta(0,T)$   the intersection of
 the domain of influence of $\Gamma$ with $\partial\Omega_0\times [0,T]$.
We assume that the domain of influence of $\Gamma$ does not intersect
$\overline{\Omega_p'}\times[0,T]$.
 
\begin{lemma}                                    \label{lma:2.1}
Suppose $\Lambda^{(1)}=\Lambda^{(2)}$ on $\Delta(0,T)$.
There exist neighborhoods
$U^{(p)}\subset \Omega^{(p)},p=1,2,
\ \overline{U^{(p)}}~\cap\partial\Omega_0=\overline{\Gamma}$
and the diffeomorphism $\varphi:U^{(1)}\rightarrow U^{(2)}$ such
that $\varphi|_\Gamma=I$  and $\|g_2^{jk}\|=\varphi\circ \|g_1^{jk}\|$.
Moreover $A_j^{(1)},\ 1\leq j\leq n,\ V^{(1)}$  and 
$\varphi\circ A_j^{(2)},1\leq j\leq n, \ \varphi\circ ~V^{(2)}$
are gauge equivalent in $U^{(1)}$, i.e. there exists 
$g(x)\in G_0(\overline{U^{(1)}}), \ g(x)=I$ on $\Gamma$ such that
(\ref{eq:1.4}) holds in $\overline{U^{(1)}}$.
\end{lemma}

The proof of Lemma \ref{lma:2.2} is the same as the proof of Lemma 2.1
in [E1].  One should replace only the inner products of the form
$\int u(x,t)\overline{v(x,t)}dxdt$ by 
$\int Tr (uv^*)dxdt$ where $v^*$ is the adjoint matrix to $v(x,t)$.
We do not assume that matrices $A_j^{(p)},V^{(p)}$ are self-adjoint
In the latter case Lemma \ref{lma:2.1} can be obtained by 
the BC-method (see[B], [KKL]).  Extend $\varphi^{-1}$ from $\overline{U_2}$
to $\overline{\Omega^{(2)}}$ in such a way that
$\varphi=I$ on $\partial\Omega_0$ and $\varphi$ is a diffeomorphism
of $\overline{\Omega^{(2)}}$ and $\overline{\tilde{\Omega}^{(2)}}=
\varphi^{-1}(\overline{\Omega^{(2)}})$.
Also extend $g(x)$ from $\overline{U_1}$ to $\overline{\tilde{\Omega}^{(2)}}$ 
so that $g(x)\in G_0(\overline{\tilde{\Omega}_2}),
 \ g=I$ on $\partial\Omega_0$.
Then we get that $\tilde{L}^{(2)}=g\circ \varphi\circ L^{(2)}=L^{(1)}$ in $U^{(1)}$.
\begin{lemma}                                  \label{lma:2.2}
Let $L^{(1)}$ and $L^{(2)}$ be the operators of the form (\ref{eq:2.4})
in $\Omega^{(p)}=\Omega_0\setminus\overline{\Omega_p'},\
p=1,2$.  Let $B\subset\Omega^{(1)}\cap\Omega^{(2)}$ 
be simply-connected,
$\partial B\cap \partial \Omega_0=\Gamma$ be 
open and connected, and  $\Omega^{(p)}\setminus\overline{B}$ be smooth.
Suppose $L^{(2)}=L^{(1)}$ in $B$ and $\Lambda^{(1)}=\Lambda^{(2)}$
on $\partial\Omega_0\times (0,T_0)$  where $\Lambda^{(p)}$  are
the D-to-N operators corresponding to $L^{(p)}, \ p=1,2.$
Then   $\tilde{\Lambda}^{(1)}=\tilde{\Lambda}^{(2)}$
where $\tilde{\Lambda}^{(p)}$ are the D-to-N  operators 
corresponding to $L^{(p)}$ in the domains  
$(\Omega^{(p)}\setminus\overline{B})\times(\delta,T_0-\delta),
\ \delta =\max_{x\in\overline{B}} d(x,\partial\Omega_0)),
\ \ d(x,\partial\Omega_0)$ is the distance in $\overline{B}$
between $x\in\overline{B}$ and $\partial\Omega_0$,
$\tilde{\Lambda}^{(p)}$ are given on 
$\partial(\Omega_0\setminus\overline{B})\times(\delta,T_0-\delta)$.
\end{lemma}

Therefore Lemma \ref{lma:2.2} reduces the inverse problem in
$\Omega^{(p)}\times(0,T_0)$ to the inverse problem in a smaller
domain 
$(\Omega^{(p)}\setminus\overline{B})\times(\delta,T_0-\delta)$.
Combining Lemmas \ref{lma:2.1} and \ref{lma:2.2} we can prove that 
for any $x^{(0)}\in \Omega^{(1)}$ there exist a 
simply-connected domain 
$B_1\subset \Omega^{(1)}$,  
 $\ x^{(0)}\in B_1$, a diffeomorphism 
$\varphi$  of $\tilde{\Omega}^{(2)}$ onto $\overline{\Omega^{(2)}}$, 
$\ \varphi=I$ on $\partial\Omega_0$,
such that
 $g\in G_0(\overline{\tilde{\Omega}^{(2)}})$ such that
$\tilde{L}^{(2)}\stackrel{\mathrm{def}}{=}
g\circ\varphi\circ L^{(2)}=L^{(1)}$ in 
$\overline{B_1}$.
To prove the global gauge equivalence and global diffeomorphism 
in the case when $\Omega^{(1)}$ is not simply-connected we shall use
some additional global quantities determined by the D-to-N
operator (c.f. [E2]).

\section{Global gauge equivalence.}
\label{section 3}
\init
 
In this section we shall prove Theorem \ref{theo:2.1}.
Fix arbitrary point $x^{(0)}\in \partial \Omega_0$.  Let $\gamma$
be a path in $\Omega$ starting at $x^{(0)}$ 
and ending at 
$x^{(1)}\in \overline{\Omega},\ \gamma(\tau)=x(\tau)$ is the
parametric equation of $\gamma,\ 0\leq\tau\leq\tau_1, \
x^{(0)}=x(0),
 \ x^{(1)}=x(\tau_1)$.  Denote by $c^{(p)}(\tau,\gamma),\ p=1,2,$
the solution of the system of differential equations
\begin{equation}                                 \label{eq:3.1}
i\frac{\partial c^{(p)}(\tau,\gamma)}{\partial\tau}=
\dot{\gamma}(\tau)\cdot A^{(p)}(x(\tau))c^{(p)}(\tau,\gamma),
\end{equation}
where 
\begin{equation}                                 \label{eq:3.2}
c^{(p)}(0,\gamma)=I_m,\ p=1,2,\ 0\leq \tau\leq \tau_1,
\end{equation}
\begin{equation}                                \label{eq:3.3}
\dot{\gamma}(\tau)=\frac{dx}{d\tau}.
\end{equation}
Denote $c^{(p)}(x^{(1)},\gamma)=c^{(p)}(\tau_1,\gamma),\ p=1,2.$
\begin{lemma}                                  \label{lma:3.1}
Suppose $A^{(1)}$ and $A^{(2)}$ are locally gauge equivalent.  Then
the
matrix $c^{(2)}(x^{(1)},\gamma)(c^{(1)}(x^{(1)},\gamma))^{-1}$ depends
only on the homotopy class of the path $\gamma$ connecting
$x^{(0)}$ and $x^{(1)}$.
\end{lemma}

\underline{Proof}
Let $\gamma_1$ and $\gamma_2$ be two homotopic paths connecting
$x^{(0)}$ and $x^{(1)}$.  Consider the path 
$\gamma_0=\gamma_1\gamma_2^{-1}$ that starts and ends at
$x^{(0)}$.  It follows from (\ref{eq:3.1})  that  
$c^{(2)}(\tau,\gamma)(c^{(1)}(\tau,\gamma))^{-1}$ satisfies
the following system of differential equations:
\begin{eqnarray}                         \label{eq:3.4}
i\frac{\partial}{\partial \tau}(c^{(2)}(\tau,\gamma)
(c^{(1)}(\tau,\gamma))^{-1})=
\dot{\gamma}(\tau)\cdot A^{(2)}(x(\tau))
(c^{(2)}(\tau,\gamma)(c^{(1)}(\tau,\gamma))^{-1})
\nonumber
\\
-(c^{(2)}(\tau,\gamma)(c^{(1)}(\tau,\gamma))^{-1})
A^{(1)}(x(\tau))\cdot \dot{\gamma}(\tau).
\end{eqnarray}
Let $b(\tau,\gamma)
=c^{(2)}(\tau,\gamma)
(c^{(1)}(\tau,\gamma))^{-1},\ b(x^{(1)},\gamma_1)
=b(\tau_1,\gamma_1),\ 
b(x^{(1)},\gamma_2)=b(\tau_2,\gamma_2)$,
where $x=x^{(p)}(\tau)$ are parametric equations of 
$\gamma^{(p)},\ 0\leq\tau\leq\tau_p,\ p=1,2$.  We have that 
$b(x^{(1)},\gamma_1)=b(x^{(1)},\gamma_2)$ iff 
$b(x^{(0)},\gamma_0)=I_m,$ where $x^{(0)}$
is the endpoint of path $\gamma_0=\gamma_1\gamma_2^{-1}$ and
$b(x^{(0)},\gamma_2)$
 is the value
at the endpoint of the solution of (\ref{eq:3.4}) along $\gamma_0$
with the initial value (\ref{eq:3.2}).  If $\gamma_0$ can be contracted 
to a point in $\overline{\Omega}$ 
there exists closed paths $\sigma_1,...,\sigma_N$ 
such that $\gamma_0=\sigma_1...\sigma_N$ and each $\sigma_j$ is 
contained in a neighborhood $U_j\subset\Omega$ where $A^{(1)}$ and
$A^{(2)}$ a gauge equivalent (see Lemma \ref{lma:2.1}).
We shall show that $b_j(\tau,\sigma_j)$ is continuous on $\sigma_j$
 where
$b_j(\tau,\sigma_j)$ is the solution of (\ref{eq:3.4}) with $\gamma$
replaced by $\sigma_j,\ \sigma_j(\tau)=x^{(j)}(\tau),
\ 0\leq \tau\leq \tau_j,$  is the parametric equation of $\sigma_j,
\ \sigma_j(0)=\sigma_j(\tau_j)$.  The continuity on $\sigma_j$ means that
$b(0,\sigma_j)=b(\tau_j,\sigma_j)$.
Since $A^{(1)}$ and $A^{(2)}$ are gauge equivalent in $U_j$ there exists 
$g_j(x)\in C^\infty(\overline{U_j})$ such that (\ref{eq:1.4}) holds in
$U_j$.  It follows from (\ref{eq:3.4}) and (\ref{eq:1.4}) that
\begin{equation}                               \label{eq:3.5}
i\frac{\partial}{\partial\tau}(b_j(\tau,\sigma_j)
g_j^{-1}(x^{(j)}(\tau))=0
\ \ \mbox{for}\ \ 0\leq\tau\leq\tau_j.
\end{equation}
Therefore $b_j(\tau,\sigma_j)
g_j^{-1}(x^{(j)}(\tau))=C$ on $\sigma_j$.
We have $b_j(0,\sigma_j)=b_j(\tau_j,\sigma_j)=C\ g_j(x^{(j)}(0)).$
Since $b_j(\tau,\sigma_j)$ is continuous on each $\sigma_j,\ 1\leq j\leq N$,
we get that $b(\tau,\gamma_0)$ is continuous on $\gamma_0$, in particular, $b(x^{(0)},0,\gamma_0)=I_m$.  Therefore  $ b(x^{(1)},\gamma_1)=
 b(x^{(1)},\gamma_2)$.
\qed

Now we shall prove that $ b(x^{(1)},\gamma_1)=
 b(x^{(1)},\gamma_2)$ for any two paths connecting $x^{(0)}$
and $x^{(1)}$.  As in the case of Lemma \ref{lma:3.1} it is enough
to prove that $b(\tau,\gamma_)$ is continuous on $\gamma_0=
\gamma_1\gamma_2^{-1}$ where $b(\tau,\gamma_0)$ is the solution of (\ref{eq:3.4})
for $\gamma_0$.

We say that $\tilde{\gamma}=\tilde{\gamma}_1,...,\tilde{\gamma}_N$ 
is a broken ray in 
$\overline{\Omega}\times[0,T_0]$
with legs $\tilde{\gamma}_j,\ 1\leq j\leq N$, if it starts at some point
$(x^{(1)},t^{(1)})\in \partial\Omega_0\times [0,T_0],\ 
x=x^{(1)}+\tau\omega,\ t=t^{(1)}+\tau$ is the parametric equation of 
$\tilde{\gamma}_1$
for $0\leq \tau\leq \tau_1$.  Then $\tilde{\gamma}$ makes $N-1$ nontangential
reflections at 
$\partial\overline{\Omega}'\times[0,T_0],$  where
$\overline{\Omega'}=\cup_{j=1}^r\overline{\Omega_j}$ and ends at
$\partial\Omega_0\times[0,T_0]$.
Denote by 
$\gamma=\gamma_1...\gamma_N$ the projection of 
$\tilde{\gamma}$ onto the $x$-plane.  
Let
$c^{(p)}(\tau,\gamma)$  be the solution of the system
\begin{equation}                                \label{eq:3.6}
i\frac{\partial}{\partial \tau}c^{(p)}(\tau,\gamma)=A^{(p)}(\gamma(\tau))\cdot
\dot{\gamma}(\tau)c^{(p)}(\tau,\gamma),
\end{equation}
\begin{equation}                                \label{eq:3.7}
c^{(p)}(0,\gamma)=I_m,\ \ \ p=1,2,
\end{equation}
where $\gamma(\tau)$ is the parametric equation of broken ray,
$0\leq \tau\leq \tau_N,\ \dot{\gamma}=\frac{d\gamma}{d\tau}$ is  the direction
of the broken ray,  $c^{(p)}(\tau,\gamma)$  is continuous on $\gamma(\tau)$.
Note that $\frac{d\gamma_j(\tau)}{d\tau}=\theta_j$ is constant on $\gamma_j,
\ \tau_{j-1}\leq \tau\leq\tau_j$.

The following lemma is the generalization of Theorem 2.1 in [E2]:
\begin{lemma}                                       \label{lma:3.2}
Let $(x^{(N)},t^{(N)})$ be the endpoint of the broken ray $\tilde{\gamma}$ :
$\tilde{\gamma}(\tau_N)=(x^{(N)},t^{(N)})$.  Denote 
$c^{(p)}(x^{(N)},\gamma)=c^{(p)}(\tau_N,\gamma),\ p=1,2$.  Then
$c^{(2)}(x^{(N)},\gamma)=c^{(1)}(x^{(N)},\gamma)$
assuming that $\Lambda^{(1)}=\Lambda^{(2)}$ on $\partial\Omega_0\times(0,T_0)$.
\end{lemma}
Assuming that Lemma \ref{lma:3.2} is proven we shall complete the proof 
of Theorem \ref{theo:2.1}

Let $\gamma$ be a broken ray starting at $x^{(1)}$ and 
ending at $x^{(N)},\ x^{(1)}\in \partial\Omega_0,
\ x^{(N)}\in\partial\Omega_0$.

  Let $c^{(p)}(\tau,\gamma)$
be the solution of (\ref{eq:3.6}), (\ref{eq:3.7}).
 Let $\alpha_1$ be a path on $\partial\Omega_0$ connecting
$x^{(0)}$ and $x^{(1)}$ and let
$\alpha_2$ be a path on $\partial\Omega_0$  connecting 
$x^{(N)}$ and $x^{(0)}$.  Therefore $\alpha=
\alpha_1\gamma\alpha_2$ is a closed path starting and ending at
$x^{(0)}$.  If $U_j\cap\partial\Omega_0\neq 0$
then the gauge $g_j=I_m$ on $\partial\Omega_0.$  Therefore
$\frac{\partial}{\partial x} g_j\cdot\stackrel{\rightarrow}{l}=0$
on $U_j\cap\partial\Omega_0$ for any vector 
$\stackrel{\rightarrow}{l}\in \R^{n}$ tangent to 
$\partial\Omega_0$.
Then (\ref{eq:1.4}) implies that $A^{(2)}\cdot \stackrel{\rightarrow}{l}
=A^{(1)}\cdot\stackrel{\rightarrow}{l} $.  It follows from (\ref{eq:3.1}),
(\ref{eq:3.2}) that $c^{(1)}(\tau,\alpha_1)=c^{(2)}(\tau,\alpha_1)$ 
and $c^{(1)}(\tau,\alpha_2)=c^{(2)}(\tau,\alpha_2)$.  Therefore $b(\tau,\alpha)$
is continuous on $\alpha=\alpha_1\gamma\alpha_2$.
We shall call $\alpha=\alpha_1\gamma\alpha_2$ an extended broken ray.

We shall assume 
for simplicity
that extended broken rays generate the homotopy group of
$\Omega$. 
Otherwise we can, as in the  end of \S 2,  construct a 
simply-connected domain 
$\Omega^{(0)}$ such that $\overline{\Omega^{(0)}}\subset\overline{\Omega},\ 
\partial\Omega^{(0)}\supset \partial\Omega_0,\ 
\partial\Omega^{(0)}\cap \overline{\Omega'}=\emptyset$  
 and
$\Omega\setminus\Omega^{(0)}$ is "thin", i.e. the volume of
$\Omega\setminus\overline{\Omega^{(0)}}$ is small.  Since 
$\Omega^{(0)}$ is homotopic to $\partial\Omega_0$
we get,  using Lemmas \ref{lma:2.1} and \ref{lma:2.2}  that potentials
$A^{(1)}$ and $A^{(2)}$ are globally gauge equivalent in 
$\overline{\Omega^{(0)}}$.  Therefore the proof of global gauge 
equivalence in $\overline{\Omega}$ can be reduced to the proof of
the global gauge equivalence in $\Omega\setminus\overline{\Omega^{(0)}}$.
It is clear that the extended broken rays in 
$\Omega\setminus\overline{\Omega^{(0)}}$ 
generate the fundamental group 
$\pi_1(\overline{\Omega}\setminus\overline{\Omega^{(0)}})$.
Note that rays without reflections also generated the fundamental group
$\pi_1(\overline{\Omega}\setminus\Omega^{(0)})$.
  Then the closed path $\gamma_0$ is homotopic  to 
$\alpha^{(1)}...\alpha^{(N_1)}$  
where $\alpha^{(j)}$ are extended broken rays.  Since $b(\tau,\alpha^{(j)})$ is
continuous on $\alpha^{(j)},\ j=1,...,N_1,$  we get that $b(\tau,\gamma_0)$
is continuous on $\gamma_0$.
It follows from Lemma \ref{lma:3.2} that 
 $b(\tau,\sigma_j)=I_m$ on $\partial\Omega_0$ and hence
$b(\tau,\gamma_0)=I_m$ on $\partial\Omega_0$.
Therefore we proved that
$c^{(2)}(x^{(1)},\gamma)(c^{(1)}(x^{(1)},\gamma))^{-1}$
does not depend on the path $\gamma$ connecting 
$x^{(0)}$ and $x^{(1)}$.  Denote $g(x^{(1)})=
c^{(2)}(x^{(1)},\gamma)(c^{(1)}(x^{(1)},\gamma))^{-1}$.
We have that $g(x)$ is a single-valued matrix on $\overline{\Omega}$,
$g(x)=I_m$ on $\partial\Omega_0$ and $g(x)$ is
nonsingular since $c^{(p)}(x,\gamma)$ are nonsingular, $p=1,2$.
We have for arbitrary $x^{(1)}$:
\begin{eqnarray}                              \label{eq:3.8}
A^{(2)}(x(\tau))\cdot\dot{\gamma}(\tau)=
i\frac{\partial c^{(2)}(\tau,\gamma)}{\partial\tau}(c^{(2)}(\tau,\gamma))^{-1}
\\
=(i\frac{\partial}{\partial x} g(x(\tau))\cdot\dot{\gamma}(\tau)c^{(1)}+ig(x(\tau))
\frac{\partial c^{(1)}}{\partial\tau})(c^{(1)}(\tau,\gamma))^{-1}g^{-1}
\nonumber
\\
=i\frac{\partial}{\partial x} g\cdot \dot{\gamma}g^{-1}+
gA^{(1)}\cdot\dot{\gamma}(\tau)g^{-1}=
(i\frac{\partial}{\partial x} g g^{-1}+gA^{(1)}g^{-1})\cdot \dot{\gamma}(\tau).
\nonumber
\end{eqnarray}

Since we can choose $\gamma(\tau)$ such that $\gamma(\tau_1)=x^{(1)}$
and $\dot{\gamma}(\tau)$ is arbitrary at $\tau=\tau_1$,  we get that
\begin{equation}                                \label{eq:3.9}
A^{(2)}(x^{(1)})=
i(\frac{\partial}{\partial x} g)g^{-1}(x^{(1)})+
g(x^{(1)})A^{(1)}(x^{(1)})g^{-1}(x^{(1)}),
\end{equation}
i.e. $A^{(2)}$ is gauge equivalent to $A^{(1)}$ in $\overline{\Omega}$.
We can change $A^{(1)}$ to $A'=gA^{(1)}g^{-1}+i\frac{\partial}{\partial x} g g^{-1},
\ V'=gV^{(1)}g^{-1}.$  Then we will have $A^{(2)}=A^{(1)}$.
Therefore applying Lemmas \ref{lma:2.1} and \ref{lma:2.2} we get that
$V^{(2)}=V'$ in $\overline{\Omega}$.  Therefore $V^{(2)}=gV^{(1)}g^{-1}$ where
$g(x)$ is the same as in (\ref{eq:3.9}).
\qed

It remains to prove Lemma \ref{lma:3.2}.  In the case when the broken ray 
$\gamma=\gamma_1...\gamma_M$ does not contain caustics points the proof
of Lemma \ref{lma:3.2} is the same as the proof of Theorem 2.1 in [E2].
We shall consider the case
when $\gamma$ has some caustics points and we shall simplify also 
the proof of
Theorem 2.1 in [E2].  
However in this paper we shall not use rays having caustics points.
Consider,  for simplicity,  the case $n=2$
and $x_*\in\gamma_M$ is the only caustics point on $\gamma$.
We also assume that the caustics point is generic (see [V]).  
Note that if $x_*$ is not generic but the broken ray $\gamma$
can be approximated by a sequence of broken rays  having 
only generic
caustics points,  then Lemma \ref{lma:3.2} 
holds for such $\gamma$ too.  This fact suggests that Lemma \ref{lma:3.2} 
is likely true for any broken ray.

Let $\chi_0(y')\in C_0^\infty(\R^2),\ y'=(y_1,y_2),\ \chi_0(y')\geq 0,\ 
\chi_0(y')=0$ for $|y'|>1,\ \chi_0(y')=1$ for $|y'|<\frac{1}{2},\ 
\int_{\R^2}\chi_0^2(y')dy'=1$.
Denote
\begin{equation}                                    \label{eq:3.10}
\chi(y')=\frac{1}{\e}\chi_0(\frac{y'}{\e}).
\end{equation}
We shall heavily use the notations of [E2, \S 2].
The difference with [E2] is that in this paper we consider the broken ray
$\tilde{\gamma}$
in $\Omega\times[0,T_0]$ and its projection on $\Omega$ will be the broken
ray $\gamma$ considered in [E2].

Let $\Pi$ be a plane in $\R^2\times \R,\ (x,t)\in\Pi$
if $x=x_0^{(0)}+y_1\omega_\perp,\ t=y_2+t^{(0)}$,  where 
$\omega_\perp\cdot \omega=0,\ x^{(0)}\not\in \Omega$ and the plane $\Pi$
does not intersect $\overline{\Omega}\times \R$.  We denote by 
$\tilde{\gamma}(y')=
\tilde{\gamma}_0(y')\tilde{\gamma}_1(y')...\tilde{\gamma}_M(y')$ 
the broken ray starting at 
$(x^{(0)}+y_1\omega_\perp,t^{(0)}+y_2)$ in the direction
$(\omega,1),\ y'=(y_1,y_2)$.   Then the equation of 
$\tilde{\gamma}_0(y')$ is
$x=x^{(0)}+y_1\omega_\perp +t_0\omega,\ 
t=t^{(0)}+y_2+t_0,\ 0\leq t_0\leq t_0(y_1),$
where $(x^{(0)}+y_1\omega_\perp+t_0(y_1)\omega, t^{(0)}+y_2+t_0(y_1))
=\tilde{P}_1$
is the point where $\tilde{\gamma}_0$ hits $\partial\Omega'\times(0,T_0)$. 
 As in [E2] we
introduce "ray coordinates" $(s_p,t_p)$ 
in the neighborhood of $\gamma_p,\ 0\leq p\leq M$.  Denote by
$D_j(x(s_j,t_j))$ the Jacobian of the change of coordinates 
$x=x^{(j)}(s_j,t_j)$.  Let $\tilde{P}_j$ be  the points of reflections of 
$\tilde{\gamma}(y')$ at $\partial\Omega',\ 1\leq j\leq M.$
Denote by $P_j$ the projection of $\tilde{P}_j$ on the $x$-plane.
Note that the time coordinate of $\tilde{P}_j$ is $t^{(j)}=t^{(0)}+y_2+
\sum_{r=0}^j t_r(y_1)$ where $t_r(y_1)$ is the distance between
 $P_r$ and $P_{r-1}$.
Note that $t=t^{(j)}+t_j$ on $\tilde{\gamma}_j,\ 0\leq t_j\leq t_j(y_1)$.

Let $L^{(p)}$ be the same as in (\ref{eq:2.1}), $p=1,2$.  We construct
a solution of $L^{(1)}u=0$ of the form (c.f. (2.1), (2.9) in [E2],
see also the earlier work [I]):
\begin{equation}                                    \label{eq:3.11}
u(x,t,\omega)=
\sum_{j=0}^{M-1}u_j(x,t,\omega)+u_{M1}+u_{M2}+u_{M3}+u^{(1)},
\end{equation}
where the principal part of $u_j$ has a form
\begin{equation}                               \label{eq:3.12}
u_{j0}=a_{j0}(x,t,\omega)e^{ik(\psi_j(x,\omega)-t)},
\end{equation}
$\psi_j(x,\omega)$ are the same as in (2.2), (2.3), (2.4) in [E2]
and
\begin{equation}                              \label{eq:3.13}
a_{j0}=(a_{j-1,0}|D_j|^{\frac{1}{2}})\left|_{P_j}\right.
\frac{1}{|D_j|^{\frac{1}{2}}}c_{j1}(x,\omega),
\end{equation}
where $c_{j1}(x,\omega)$ is the solution of the system
\begin{equation}                              \label{eq:3.14}
i\theta_j\cdot\nabla c_{j1}= (A^{(1)}\cdot \theta_j)c_{j1},
\ \ \ t_{j-1}(y_1)\leq t_j\leq t_j(y_1),\ \ \ c_{j1}\left|_{P_j}\right.
=I_m,
\end{equation}
$0\leq j\leq M,\ \nabla
=\frac{\partial}{\partial x},
\ \theta_j$ is the direction of $\gamma_j,
\ \theta_0=\omega$.

We shall assume that 
\begin{equation}                         \label{eq:3.15}
u_0(x,t,\omega)=\chi(y')\alpha_0
\end{equation}
on the plane $\Pi$, i.e. when $t_0=0$ and $x=x^{(0)}+y_1\omega_\perp,\ 
t=t^{(0)}+y_2$.  Here $\alpha_0$ is an arbitrary constant matrix.

Let $(x_*,t_*) \in \tilde{\gamma}_M$ be such that $x_*$ is the caustics point
in the $x$-plane.
Note that $u_{M1}$ has the same form as $u_{M-1}$ for $t<t_*-C\e$,
where $\e$ is the same as in (\ref{eq:3.10}),  solution $u_{M2}$ is defined
in a $C\e$-neighborhood $U_{\e}$ of $(x_*,t_*)$.  We will not write 
the explicit form of $u_{M2}$ (see, for example,  [V]) since we will
only need an estimate
\begin{equation}                               \label{eq:3.16}
|u_{M2}|\leq\frac{Ck^{\frac{1}{6}}}
{1+k^{\frac{1}{6}}d^{\frac{1}{4}}(x)},
\end{equation}
where $d(x)$ is the distance from $x\in U_{0,\e}$ to the caustics curve.
Such estimate holds in the generic case (see  [V]).
Moreover,
\[
|\nabla u_{M2}|\leq\frac{Ck^{\frac{7}{6}}}
{1+k^{\frac{1}{6}}d^{\frac{1}{4}}(x)}.
\]
Finally,  $u_{M3}$ is defined for $t>t_*+C\e$ and it has
the same form as $u_{M1}$.  The main difference is that the amplitude of
$u_{M3}$ has an extra factor $e^{i\beta}$ where $\beta$ is real.
The construction and the estimate of $u^{(1)}$ in (\ref{eq:3.11}) is similar to
[E2, Lemma 2.1] with the simplification that we consider the hyperbolic
initial-boundary value problem
with the zero initial conditions when $t=0$ and zero boundary conditions
on $\partial\Omega\times(0,T_0)$
 instead of (2.9) in [E2].  Since we assumed
that $T_0$ is large enough we get that the endpoint of 
$\tilde{\gamma}_M$ belongs to
$\partial\Omega_0\times (0,T_0)$.

We construct a solution $v(x,t,\omega)$ of $L_2^*v=0$ similar to (\ref{eq:3.11})
with the same initial data as (\ref{eq:3.15}) for $v_0$ with $\alpha_0$
replaced 
by $\beta_0$ where $\beta_0$ is an aritrary constant matrix and with the same
phase function $\psi_j(x,\omega),\ 0\leq j\leq M$,  as in (\ref{eq:3.12}):
We have
\begin{equation}                                  \label{eq:3.17}
v=\sum_{j=0}^{M-1} v_j(x,t,\omega)+v_{M1}+v_{M2}+v_{M3}+
v^{(1)}(x,t,\omega),
\end{equation}
where the principal term of $v_j$ has the following form:
\begin{equation}                               \label{eq:3.18}
v_{j0}=b_{j0}(x,t,\omega)e^{ik(\psi_j(x,\omega)-t}),
\end{equation}
where  $b_{j0}$ are the same as $a_{j0}$ with 
$c_{j1}(x,\omega)$ replaced by $c_{*,j}(x,\omega)$ where $c_{*j}$ 
is the solution of the system
\begin{equation}                                \label{eq:3.19}
i\theta_j\cdot \nabla c_{*j}=(\theta_j\cdot (A^{(2)})^*)c_{*j}.
\end{equation}
Taking the adjoint of (\ref{eq:3.19}) we get
\begin{equation}                                \label{eq:3.20}
-i\theta_j\cdot \nabla c_{*j}^*=c_{*j}^*(\theta_j\cdot A^{(2)}).
\end{equation}
Denote
\begin{equation}                               \label{eq:3.21}
c_{j2}=(c_{*j}^*)^{-1}.
\end{equation}
Then (\ref{eq:3.20}) implies that 
\begin{equation}                                \label{eq:3.22}
i\theta_j\cdot \nabla c_{j2}=(\theta_j\cdot A^{(2)})c_{j2}.
\end{equation}
We assume that $v^{(1)}$ satisfies zero initial conditions when
$t=T,\ x\in\Omega$, and zero boundary condiions on $\partial\Omega\times(0,T_0)$.
Substitute (\ref{eq:3.11}) instead of $u^{(1)}$ and (\ref{eq:3.17})
instead of $v^{(2)}$ in the Green's formula.
Dividing by $2k$ and passing to the limit when $k\rightarrow\infty$
we obtain  (c.f. [E2]):
\begin{equation}                             \label{eq:3.23}
0=\sum_{j=0}^{M-1}\int_\Omega\int_0^T((A^{(1)}-A^{(2)})\cdot
\nabla(\psi_j-t)a_{j0},b_{j0})dxdt +I_M,
\end{equation}
where $I_M$ is the integral over a neighborhood of $\gamma_M$.
We make a series of changes of variables as in (2.43) in [E2].

Note that the Jacobian $D_M(x^{(M)}(s_M,t_M))$ vanishes on the caustics set
and therefore $D_M^{-1}$ has a singularity there.  However when we make 
changes of variables this singularity in $u_{M1},v_{M1}$ and in
$u_{M3},v_{M3}$ cancels.  Note also that the estimate (\ref{eq:3.16})
implies that the integral over the neighborhood $U_\e$ is $O(\sqrt{\e})$.
Therefore taking into account that $\alpha_0$  and $\beta_0$ are
arbitrary matrices we get
\begin{equation}                              \label{eq:3.24}
\sum_{j=0}^M\int_{\R^2}\int_{\tilde{\gamma}_j(y')}\chi^2(y')
c_{j2}^{-1}(A^{(1)}-A^{(2)})
\cdot\theta_j c_{j1}dt_jdy'+O(\sqrt{\e})=0,
\end{equation}
where $\tilde{\gamma}(y')$ is the broken ray starting at 
$(x^{(0)}+y_1\omega_\perp,t^{(0)}+y_2)$
and we use in (\ref{eq:3.24}) that 
$c_{*j}^*=c_{j2}^{-1}$ 
(see (\ref{eq:3.21})).  Note that
 \begin{equation}                               \label{eq:3.25}
c_{j2}^{-1}(A^{(1)}-A^{(2)})\cdot \theta_jc_{j1}=i\theta_j\cdot 
\nabla(c_{j2}^{-1}c_{j1}),
\end{equation}
since $-c_{j2}^{-1}(A^{(2)}\cdot \theta_j)=i\theta_j\cdot 
\nabla c_{j2}^{-1}$.
After changes of variables $c_{j1}$ and $c_{j2}$ in 
(\ref{eq:3.24}) satisfy the differential equations 
(\ref{eq:3.14}), (\ref{eq:3.22}) but the initial conditions are
different : 
\begin{equation}                               \label{eq:3.26}
c_{ji}\left|_{P_j}=c_{j-1,i}\right|_{P_j},
\ \ \  1\leq j\leq M,\ \ i=1,2.  
\end{equation}
We kept the same notation for the simplicity. 
Taking the limit in (\ref{eq:3.24}) when $\e\rightarrow 0$ we get
\[
\sum_{j=0}^M\int_{\gamma_j}\theta_j\cdot\nabla(c_{j2}^{-1}c_{j1})dt_j=
\sum_{j=0}^M\left[(c_{j2}^{-1}c_{j1})\left|_{P_j}\right.-
(c_{j2}^{-1}c_{j1})\left|_{P_{j-1}}\right.\right]=0.
\]
Since $c_{0i}\left|_{P_0}\right.=I_m$ and (\ref{eq:3.26}) holds
we get that 
$c_{M2}^{-1}c_{M1}\left|_{P_{M}}\right.=I_m$,  i.e.
$c_{M1}\left|_{P_M}\right.=c_{M2}\left|_{P_M}\right.$.
Lemma \ref{lma:3.2} is proven.
\qed

\section{Global diffeomorphism.}
\label{section 4}
\init

Let $L^{(1)}u^{(1)}=0$ and $L^{(2)}u^{(2)}=0$ be equations of the form 
(\ref{eq:2.4}) in domains $\Omega^{(p)}=\Omega_0\setminus\Omega_p'$,  where
$\Omega_p'=\cup_{j=1}^r\Omega_{jp},\ p=1,2$.
We assume that the initial conditions (\ref{eq:2.2}) in $\Omega^{(p)},\ p=1,2$
and the boundary
conditions (\ref{eq:2.3}) with $\Omega_j$ replaced by $\Omega_{jp},\ p=1,2,$
are satisfied.

\begin{theorem}                                   \label{theo:4.1}
Suppose $\Lambda^{(1)}=\Lambda^{(2)}$ on $\partial\Omega_0$,
where $\Lambda^{(p)}$ are the D-to-N operators corresponding to 
$L^{(p)},p=1,2.$
Suppose 
\[
T_0>2\min_{p}\max_{x\in \overline{\Omega^{(p)}}}d_p(x,\partial\Omega_0)
\]
 where $d_p$ is the distance with respect to the metric tensor $\|g_p^{jk}\|^{-1}$.  
Then there exists a diffeomorphism $\varphi$ of $\overline{\Omega^{(1)}}$ onto  
$\overline{\Omega^{(2)}}$
such that $\varphi=I$ on $\partial\Omega_0$ and 
$\|g_2^{jk}\|=\varphi\circ\|g_1^{jk}\|$.
\end{theorem}
We shall sketch the proof of Theorem \ref{theo:4.1} assuming for
the simplicity that $m=1,\ A_j^{(p)}\equiv 0,\ 1\leq j\leq r,\ p=1,2,$
and $T_0=\infty$.  By using Lemmas \ref{lma:2.1} and \ref{lma:2.2}
we can get a 
simply-connected
domain $\Omega^{(0)}\subset\Omega^{(1)}$ such that 
$\Omega^{(1)}\setminus\Omega^{(0)}$ 
has a small volume.  Moreover there exists 
a diffeomorphism $\tilde{\varphi}$ of $\Omega^{(2)}$ onto
$\tilde{\Omega}^{(2)}=\tilde{\varphi}^{-1}(\Omega^{(2)}),\ 
\tilde{\varphi}=I$ on $\partial\Omega_0$  such that 
$\tilde{L}^{(2)}\stackrel{\mathrm{def}}{\equiv}\tilde{\varphi}\circ L_2$ is
equal to $L^{(1)}$ in $\Omega^{(0)}$.  Note that  
$\Omega^{(0)}\subset \Omega^{(1)}\cap \tilde{\Omega}^{(2)}$.  We also get 
from Lemma \ref{lma:2.2} that $\Lambda^{(1)}=\tilde{\Lambda}^{(2)}$ on 
$\partial\Omega^{(0)}\setminus  \partial\Omega_0$ 
where $\tilde{\Lambda}^{(2)}$ is the D-to-N
operator corresponding to $\tilde{L}^{(2)}$.  Since
$\Omega^{(1)}\setminus\Omega^{(0)}$ is thin, there is an open subset
$\Gamma_1$ of $\partial\Omega^{(0)}$ such that 
the endpoints of geodesics corresponding to $L^{(1)}$ in
$\Omega^{(1)}\setminus\Omega^{(0)}$, orthogonal to $\Gamma_1$,
form an open
subset $\Gamma_2\subset\partial\Omega^{(0)}$.  Denote 
by $D_1\subset\Omega^{(1)}\setminus\overline{\Omega^{(0)}}$ 
the union of these geodesics.  It follows from the proof
of Lemma \ref{lma:2.1} (see [E1]) that $\Lambda^{(1)}$ on $\Gamma_1$
uniquely determines the metric tensor $\|g_1^{jk}\|^{-1}$ in the
semi-geodesic coordinates in $D_1$.  Denote by $\psi_1$
   the map  
of $D_1$ on $\tilde{D}_1=\psi_1(D_1)$ such that 
$\psi_1(x)$ are
the semi-geodesic coordinates in $\tilde{D}_1$. 
Analogously let $D_2$ be the union of all geodesics of $\tilde{L}^{(2)}$ 
orthogonal to $\Gamma_1$ and let $\Gamma_2'\subset \partial\Omega^{(0)}$
be the set of its endpoints.  Denote by $\psi_2(x)$ 
the semi-geodesic coordinates for $\tilde{L}^{(2)}$ and let
$\tilde{D}_2=\psi_2(D_2)$.
 By Lemma \ref{lma:2.1} 
$\psi_1\circ L^{(1)}=\psi_2\circ \tilde{L}^{(2)}$
in $\tilde{D}_1\cap\tilde{D}_2$.
It follows from Lemma \ref{lma:2.1}
 that 
$\psi_j=I$ on $\Gamma_1$.  
  Note that $\Omega^{(2)}\setminus\Omega^{(0)}$ 
coincide with $\Omega^{(1)}\setminus\Omega^{(0)}$ near $\Gamma_1$.
\begin{lemma}                             \label{lma:4.1}
The following equalities hold:
$\tilde{D}_1=\tilde{D}_2$ and
$\psi=\psi_2^{-1}\psi_1=I \ \mbox{on}\ \  \Gamma_2$.
\end{lemma}
\underline{Proof:}
Since we assume that $T_0=\infty$ we can switch to the inverse problem 
for the equations of the form (\ref{eq:1.1}).  Choose parameter $k\in\C$
such that the boundary value problem of the form
(\ref{eq:1.1}), (\ref{eq:1.2}), (\ref{eq:1.3})  has a unique solution 
$u_p$
for any 
$f\in H_{\frac{1}{2}}(\partial\Omega^{(0)}\setminus\partial\overline{\Omega_0})$
 where $f$ is the same for $p=1$ and $p=2$.
Choose f nonsmooth.  
Denote $\tilde{\Gamma}_2=\psi_1^{-1}(\Gamma_2),\
\tilde{\Gamma}_2'=\psi_2^{-1}(\Gamma_2')$.
It follows from the unique continuation theorem that 
$\psi_1\circ u_1=\psi_2\circ u_2$ 
in $\tilde{D}_2\cap\tilde{D}_1$
since the Cauchy data of $u_1$ and $u_2$ 
coincide on $\Gamma_1$.  Here $L^{(p)}u_p=0,\ p=1,2$.  If
$\tilde{\Gamma}_2\neq\tilde{\Gamma}_2'$ 
we get a contradiction since 
$\psi_2\circ u_1$ is $C^\infty$ outside $\tilde{\Gamma}_2$ and 
$\psi_1\circ u_2$ is $C^\infty$ outside $\tilde{\Gamma}_2'$.  
Therefore $\tilde{D}_1=\tilde{D}_2$  and $\psi=\psi_2^{-1}\psi_1$ is 
a diffeomorphism of $D_1$ onto $D_2$.  
Since $\tilde{\Gamma}_2=\tilde{\Gamma}_2'$ we have 
$\Gamma_2'=\psi(\Gamma_2)$.  Since $\psi_1\circ u_1=\psi_2\circ u_2$ in
$\tilde{D}_1=\tilde{D}_2$ and $u_1=u_2=f(x)$ on $\Gamma_2$ we get
$f(x)=f(\psi(x))$ on $\Gamma_2$.
Since $f$ is arbitrary this implies that $\psi=I$ on $\Gamma_2.$
\qed

Therefore 
$\psi=I$ on $\partial D_1\cap \partial \Omega^{(0)}$.
Define $\varphi^{(1)}=\tilde{\varphi}$ on $\Omega^{(0)},\ \varphi^{(1)}=
\psi\circ \tilde{\varphi}$ on $D_1$.  We get that 
$\varphi^{(1)}\circ L^{(2)}=L^{(1)}$
in $\overline{\Omega^{(0)}}\cup\overline{D_1}$.

Applying Lemma \ref{lma:2.2} to 
$\Omega^{(p)}\setminus(\overline{\Omega^{(0)}}\cup\overline{D_1})$ and using again
Lemmas \ref{lma:4.1}, \ref{lma:2.1} and \ref{lma:2.2} we prove
Theorem \ref{theo:4.1}.

\underline{Remark 4.1} (c.f. [E1]).
We shall show now that the obstacles can be recovered up to 
the diffeomorphism.  Let $\gamma_0$ be an open subset of 
$\partial\Omega^{(0)}$ close to the obstacle $\Omega_1'$.
Denote by $\Delta_1$ the union of all geodesics in $\Omega^{(1)}$
orthogonal to $\gamma_0$ and ending on $\Omega_1'$.  Denote by
$\gamma_1$ the intersection of $\overline{\Delta_1}$ and 
$\overline{\Omega_1'}$.
Introduce semi-geodesic coordinates for $L^{(1)}$ in $\Delta_1$. 
Let $\varphi_1$ be the change of variables to the semi-geodesic
coordinates and let $\tilde{\Delta}_1=\varphi_1(\Delta_1)$.
Let $\varphi_2$ be the change of variables to the semi-geodesic
coordinates for $\tilde{L}^{(2)}$ in 
$\Delta_2$  where $\Delta_2$ is the union of all geodesics of
$\tilde{L}^{(2)}$ orthogonal to $\gamma_0$ and ending on 
$\Omega_2'$.  Let $\gamma_2=
\overline{\Delta_2}\cap\partial\Omega_2',\ \tilde{\gamma}_2=
\varphi_2(\gamma_2),\ 
\tilde{\Delta}_2=\varphi_2(\Delta_2)$.
Let $L^{(1)}u_1=0$ be a geometric optics solution in
$\Omega^{(1)}\setminus\Omega^{(0)}$ similar to constructed in
\S 3 that starts on $\gamma_0$, reflects at $\partial\Omega_i'$ and
leaves $\Omega^{(1)}\setminus\Omega^{(0)}$ again on $\gamma_0$.
Let $u_2$ be the solution of $\tilde{L}^{(2)}u_2=0$ in
$\Omega^{(2)}\setminus\Omega^{(0)}$ having the same boundary data as 
$u_1$.  Since $\varphi_1\circ L^{(1)}=\varphi_2\circ\tilde{L}^{(2)}$ in
$\tilde{\Delta}_1\cap\tilde{\Delta}_2$ 
and since $\varphi_1\circ u_1$ and $\varphi_2\circ u_2$
have the same Cauchy data on $\gamma_0$ we get by the uniqueness
continuation theorem that $\varphi_1\circ u_1=\varphi_2\circ u_2$
in $\tilde{\Delta}_1\cap\tilde{\Delta}_2$.
If
$\tilde{\gamma}_1\neq\tilde{\gamma}_2$  then 
we can find $u_1$ such that
$\varphi_1\circ u_1$
and $\varphi_2\circ u_2$ will have different point of reflection and this 
will contradict that  
$\varphi_1\circ u_1=\varphi_2\circ u_2$
in $\tilde{\Delta}_1\cap\tilde{\Delta}_2$. 
  Since $\tilde{\gamma}_1=\tilde{\gamma}_2$ we
get that $\varphi(\gamma_1)=\gamma_2\subset\partial\Omega_2'$
and $\varphi(\Delta_1)=\Delta_2$  where $\varphi=\varphi_2^{-1}\varphi_1$.
\qed

\underline{Remark 4.2}
Note that Lemma \ref{lma:4.1} allows to consider the inverse problems 
in multi-connected domains $\Omega$ with the D-to-N operator given on
a not connected part $\Gamma_0$ of $\partial\Omega$.

\end{document}